\documentclass[12pt]{article}
\usepackage[english]{babel}
\usepackage{latexsym}
\usepackage{amssymb,amsthm,amsmath}
\usepackage{graphicx}
\title{\bf A computer based classification of caps in $PG(4,2)$}
\author{Daniele Bartoli, Stefano Marcugini, Fernanda Pambianco}

{\theoremstyle{definition}
\newtheorem*{definition*}{Definition}
}

\newtheorem*{proposition*}{Proposition}
\newtheorem*{corollary*}{Corollary}
\newtheorem*{lemma*}{Lemma}
\date{}
\begin{document}

\maketitle \vspace*{-15mm} \noindentÄ

 \vspace*{2cm}
\begin{abstract}
\noindent In this paper we present the complete classification of
caps in $PG(4,2)$. These results have been obtained using a computer
based exhaustive search that exploits projective equivalence.
\end{abstract}

\section{Introduction}
In the projective space $PG(r,q)$ over the Galois Field $GF(q)$, a $n$-cap is a set of $n$ points
 no 3 of which are collinear. A $n$-cap is called \emph{complete} if it is not contained in a
 $(n+1)$-cap. For  a detailed description of the most important properties of these geometric
 structures, we refer the reader to \cite{Hir}.
 In the last decades the problem of determining the spectrum of the
 sizes of complete caps has been the subject of a lot of researches.
 For a survey see \cite{18secondo}.\\
 In \cite{Hir2} (see Th.  $18.2.1$) is presented the classification of complete caps in $PG(3,2)$ and it is also possible to understand the classification of the incomplete caps.  The following table shows the number of non equivalent complete and incomplete caps in $PG(3,2)$.
 \newpage
 \begin{table}[h]
\caption{Number and type of non equivalent examples}
\begin{center}
\begin{tabular}{|c|c|c|}
\hline
$|\mathcal{K}|$&\# COMPLETE&\# INCOMPLETE\\
&CAPS&CAPS\\
\hline
5&1&1\\
6&0&1\\
7&0&1\\
8&1&0\\
\hline
\end{tabular}
\end{center}
\end{table}

In this work we search for the classification of complete and
incomplete caps in $PG(4,2)$, using an exhaustive search algorithm. In Section \ref{alg}
the algorithm utilized is illustrated;
 in Section \ref{risultati} the complete list of non equivalent complete and incomplete
caps is presented.

\section{The searching algorithm }\label{alg}
In this section the algorithm utilized is presented. Our goal is to obtain the classification of complete and incomplete caps in $PG(4,2)$. It is not restrictive to suppose that a cap in $PG(4,2)$ contains this five points:
\begin{displaymath}
\mathcal{R}=\{(1:0:0:0:0);(0:1:0:0:0);(0:0:1:0:0);(0:0:0:1:0);(0:0:0:0:1)\}.
\end{displaymath}
Then we define the set \emph{Cand} of all the points lying no 2-secant of $\mathcal{R}$. We introduce in \emph{Cand} the following equivalence relationship:
\begin{displaymath}
P \sim Q \iff \mathcal{C} \cup \{P\} \cong \mathcal{C} \cup \{Q\},
\end{displaymath}

\noindent where $\cong$ means that the two sets are projectively equivalent. This relationship spreads the candidates in equivalent classes $\mathcal{C}_{1},\ldots,\mathcal{C}_{k}$.\\
The choice of the next point to add to the building cap can be made only among the representatives of the equivalent classes, in fact two caps one containing $\mathcal{C} \cup \{P\}$  and the other one $\mathcal{C} \cup \{Q\}$, with $P$ and $Q$ in $\mathcal{C}_{\overline{i}}$, are equivalent by definition of orbit.\\
Suppose now that we have construct all the caps containing $\mathcal{C} \cup \{P_{i}\}$, with $i \leq \overline{i}$. Considering the caps containing $\mathcal{C} \cup \{P_{j}\}$ with $\overline{i}<j$, all the points of the classes $\mathcal{C}_{k}$ with $k<\overline{i}$ can be avoided. In fact a cap containing $\mathcal{C} \cup \{P_{j}\} \cup \{\overline{P}_{k}\}$, with $\overline{P}_{k} \in \mathcal{C}_{k}$ and $k<\overline{i}$, is projectively equivalent to a cap containing $\mathcal{C} \cup \{P_{k}\} \cup \{P_{j}\}$, already studied.\\
When we add a new point to the cap, we can divide all the remaining candidates in equivalence classes, as above. Two points  $P$ and $Q$ are in relationship with the $j$-th class $\mathcal{C}_{j}$, i.e. $P \sim_{j} Q$, if $\mathcal{C} \cup \{P_{j}\} \cup \{P\}$ and $\mathcal{C} \cup \{P_{j}\} \cup \{Q\}$ are projectively equivalent. \\
At the $m$-th step of the extension process if the cap $\mathcal{C}\cup \{P_{i_{m}}\} \cup \ldots \cup \{P^{i_{1}\ldots i_{m-1}}_{i_{m}}\} \cup\{P\}$  is projectively equivalent to the cap $\mathcal{C}\cup \{P_{i_{m}}\} \cup \ldots \cup \{P^{i_{1}\ldots i_{m-1}}_{i_{m}}\} \cup\{Q\}$ with $P^{i_{1}\ldots i_{r}}_{s} \in \mathcal{C}^{i_{1}\ldots i_{r}}_{s} $, then $P$ and $Q$ are in relationship ($P \sim_{i_{1}\ldots i_{m}} Q$) and they belong to the same class $\mathcal{C}^{i_{1}\ldots i_{m}}_{m+1} $.\\
Iterating the process we can build a tree similar to the following:
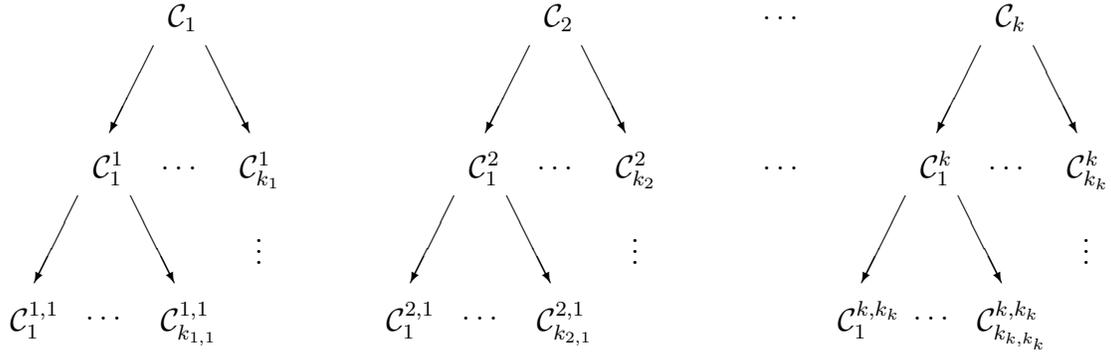
\begin{center}

\setlength{\unitlength}{1cm}
\begin{picture}(15,6)(0,0)


\put(2.35,4.9){$\mathcal{C}_{1}$}

\put(2.15,4.65){\vector(-1,-2){0.58}}
\put(2.85,4.65){\vector(1,-2){0.58}}

\put(1.35,2.9){$\mathcal{C}_{1}^{1}$}
\put(2.25,3){\ldots}

\put(3.3,2.9){$\mathcal{C}_{k_{1}}^{1}$}

\put(1.15,2.65){\vector(-1,-2){0.58}}
\put(1.85,2.65){\vector(1,-2){0.58}}

\put(3.5,1.75){\vdots}

\put(0.25,0.85){$\mathcal{C}_{1}^{1,1}$}
\put(1.25,1){\ldots}

\put(2.25,0.85){$\mathcal{C}_{k_{1,1}}^{1,1}$}

\put(10.25,5){\ldots}
\put(10.25,3){\ldots}


\put(7.35,4.9){$\mathcal{C}_{2}$}

\put(7.15,4.65){\vector(-1,-2){0.58}}
\put(7.85,4.65){\vector(1,-2){0.58}}

\put(6.35,2.9){$\mathcal{C}_{1}^{2}$}
\put(7.25,3){\ldots}

\put(8.30,2.9){$\mathcal{C}_{k_{2}}^{2}$}

\put(6.15,2.65){\vector(-1,-2){0.58}}
\put(6.85,2.65){\vector(1,-2){0.58}}

\put(8.5,1.75){\vdots}

\put(5.25,0.85){$\mathcal{C}_{1}^{2,1}$}
\put(6.25,1){\ldots}

\put(7.25,0.85){$\mathcal{C}_{k_{2,1}}^{2,1}$}


\put(13.35,4.9){$\mathcal{C}_{k}$}

\put(13.15,4.65){\vector(-1,-2){0.58}}
\put(13.85,4.65){\vector(1,-2){0.58}}

\put(12.35,2.9){$\mathcal{C}_{1}^{k}$}
\put(13.25,3){\ldots}

\put(14.3,2.9){$\mathcal{C}_{k_{k}}^{k}$}

\put(12.15,2.65){\vector(-1,-2){0.58}}
\put(12.85,2.65){\vector(1,-2){0.58}}

\put(14.5,1.75){\vdots}

\put(11.25,0.85){$\mathcal{C}_{1}^{k,k_{k}}$}
\put(12.25,1){\ldots}

\put(13.1,0.85){$\mathcal{C}_{k_{k,k_{k}}}^{k,k_{k}}$}

\end{picture}
\end{center}
The tree is important to restrict the number of candidates in the extension process. Suppose that we have generated a $n$-cap containing the cap $\mathcal{C}\cup \{P_{i_{1}}\} \cup \{P_{i_{2}}^{i_{1}}\} \cup \ldots \cup \{P^{i_{1}\ldots i_{m-1}}_{i_{m}}\} \cup\{P\}$, after having generated $n$-caps containing $\mathcal{C}\cup \{P_{j}\}$ with $j<i_{1}$, $\mathcal{C}\cup \{P_{i_{1}}\}\cup \{P^{i_{1}}_{j}\}$ with $j<i_{2}$,\ldots, $\mathcal{C}\cup \{P_{i_{1}}\}\cup \{P^{i_{1}}_{i_{2}}\}\cup \ldots \cup \{P^{i_{1}\ldots i_{m-1}}_{j}\}$ with $j<i_{m}$, with $P^{i_{1}\ldots i_{r}}_{s} \in \mathcal{C}^{i_{1}\ldots i_{r}}_{s} $. Then the points belonging to $\mathcal{C}_{1} \cup\ldots \cup \mathcal{C}_{i_{1}-1} \cup\mathcal{C}^{i_{1}}_{1} \cup\ldots \cup \mathcal{C}^{i_{1}}_{i_{2}-1} \cup \ldots \cup \mathcal{C}^{i_{1}\ldots i_{m-1}}_{1} \cup\ldots \cup \mathcal{C}^{i_{1}\ldots i_{m-1}}_{i_{m}-1}$ can be avoided, because a cap containing one of them is equivalent to one already found. For example a $n$-cap containing $\mathcal{C}\cup \{P_{i_{1}}\} \cup \ldots \cup \{P^{i_{1}\ldots i_{m-1}}_{i_{m}}\} \cup\{P\}\cup \{Q\}$ with $Q \in \mathcal{C}_{h}$ for some $h<i_{1}$ is equivalent to a $n$-cap containing $\mathcal{C} \cup \{P_{h}\}$, which is already found. \\

\section{Results}\label{risultati}
In this Section all non equivalent caps, complete and incomplete, in
$PG(4,2)$ are presented. 
\subsection{Non-equivalent caps $\mathcal{K}$ in $PG(4,2)$}
This table shows the number and the type of the non equivalent
examples of all the caps.
\begin{table}[h]
\caption{Number and type of non equivalent examples}
\begin{center}
\begin{tabular}{|c|c|c|}
\hline
$|\mathcal{K}|$&\# COMPLETE&\# INCOMPLETE\\
&CAPS&CAPS\\
\hline
6&0&3\\
7&0&3\\
8&0&4\\
9&1&4\\
10&1&3\\
11&0&2\\
12&0&2\\
13&0&1\\
14&0&1\\
15&0&1\\
16&1&0\\
\hline
\end{tabular}
\end{center}
\end{table}
\newpage

All the non equivalent examples of the caps described in the table above are presented. In particular for each cap we compute the stabilizer and the weight enumerators of the associated code.
\begin{center}
\tabcolsep=0.85 mm
\begin{tabular}{|c|c|c|c|c|}
\hline
Size&Type& Cap&Weight enumerators&$G$\\
\hline
\hline
$6$&
Inc.&
{\scriptsize
\tabcolsep=0.85 mm
\begin{tabular}{cccccc}
 1& 0& 0& 0& 0& 1\\
 0& 1& 0& 0& 0& 1\\
 0& 0& 1& 0& 0& 1\\
 0& 0& 0& 1& 0& 1\\
 0& 0& 0& 0& 1& 1\\
\end{tabular}
}
&$2^{15}4^{15}6^{1}$&$G_{720}$\\
\hline
$6$&
Inc.&
{\scriptsize
\tabcolsep=0.85 mm
\begin{tabular}{cccccc}
 1& 0& 0& 0& 1& 0\\
 0& 1& 0& 0& 0& 0\\
 0& 0& 1& 0& 1& 0\\
 0& 0& 0& 1& 1& 0\\
 0& 0& 0& 0& 1& 1\\
\end{tabular}
}
&$1^{1}2^{10}3^{10}4^{5}5^{5}$&$G_{120}$\\
\hline
$6$&
Inc.&
{\scriptsize
\tabcolsep=0.85 mm
\begin{tabular}{cccccc}
 1& 0& 0& 0& 0& 1\\
 0& 1& 0& 0& 0& 0\\
 0& 0& 1& 0& 0& 1\\
 0& 0& 0& 1& 0& 1\\
 0& 0& 0& 0& 1& 0\\
\end{tabular}
}
&$1^{2}2^{7}3^{12}4^{7}5^{2}6^{1}$&$G_{48}$\\
\hline
\hline
$7$&
Inc.&
{\scriptsize
\tabcolsep=0.85 mm
\begin{tabular}{ccccccc}
 1& 0& 0& 0& 0& 1& 0\\
 0& 1& 0& 0& 0& 0& 1\\
 0& 0& 1& 0& 0& 1& 0\\
 0& 0& 0& 1& 0& 1& 1\\
 0& 0& 0& 0& 1& 0& 1\\
\end{tabular}
}
&$2^{6}3^{9}4^{9}5^{6}7^{1}$&$G_{72}$\\
\hline
$7$&
Inc.&
{\scriptsize
\tabcolsep=0.85 mm
\begin{tabular}{ccccccc}
 1& 0& 0& 0& 0& 1& 1\\
 0& 1& 0& 0& 0& 1& 0\\
 0& 0& 1& 0& 0& 0& 1\\
 0& 0& 0& 1& 0& 1& 1\\
 0& 0& 0& 0& 1& 1& 0\\
\end{tabular}
}
&$2^{5}3^{12}4^{7}5^{4}6^{3}$&$G_{48}$\\
\hline
$7$&
Inc.&
{\scriptsize
\tabcolsep=0.85 mm
\begin{tabular}{ccccccc}
 1& 1& 0& 0& 0& 0& 1\\
 0& 1& 1& 0& 0& 0& 0\\
 0& 1& 0& 1& 0& 0& 1\\
 0& 0& 0& 0& 1& 0& 1\\
 0& 0& 0& 0& 0& 1& 0\\
\end{tabular}
}
&$1^{1}2^{3}3^{11}4^{11}5^{3}6^{1}7^{1}$&$G_{48}$\\
\hline
\hline
$8$&
Inc.&
{\scriptsize
\tabcolsep=0.85 mm
\begin{tabular}{cccccccc}
 1& 1& 1& 0& 0& 0& 0& 1\\
 0& 1& 0& 1& 0& 0& 0& 0\\
 0& 1& 1& 0& 1& 0& 0& 1\\
 0& 0& 0& 0& 0& 1& 0& 1\\
 0& 0& 1& 0& 0& 0& 1& 0\\
\end{tabular}
}
&$2^{4}4^{22}6^{4}8^{1}$&$G_{384}$\\
\hline
$8$&
Inc.&
{\scriptsize
\tabcolsep=0.85 mm
\begin{tabular}{cccccccc}
 1& 1& 0& 0& 0& 0& 0& 1\\
 0& 1& 1& 1& 0& 0& 0& 0\\
 0& 1& 0& 1& 1& 0& 0& 1\\
 0& 0& 0& 1& 0& 1& 0& 1\\
 0& 0& 0& 0& 0& 0& 1& 0\\
\end{tabular}
}
&$1^{1}3^{7}4^{14}5^{7}7^{1}8^{1}$&$G_{168}$\\
\hline
$8$&
Inc.&
{\scriptsize
\tabcolsep=0.85 mm
\begin{tabular}{cccccccc}
 1& 1& 0& 0& 0& 0& 1& 1\\
 0& 1& 1& 0& 0& 0& 1& 0\\
 0& 1& 0& 1& 0& 0& 0& 1\\
 0& 0& 0& 0& 1& 0& 1& 1\\
 0& 0& 0& 0& 0& 1& 1& 0\\
\end{tabular}
}
&$2^{1}3^{10}4^{11}5^{4}6^{3}7^{2}$&$G_{48}$\\
\hline
$8$&
Inc.&
{\scriptsize
\tabcolsep=0.85 mm
\begin{tabular}{cccccccc}
 1& 1& 0& 0& 0& 0& 1& 0\\
 0& 1& 1& 0& 0& 0& 0& 1\\
 0& 1& 0& 1& 0& 0& 1& 0\\
 0& 0& 0& 0& 1& 0& 1& 1\\
 0& 0& 0& 0& 0& 1& 0& 1\\
\end{tabular}
}
&$2^{2}3^{8}4^{10}5^{8}6^{2}8^{1}$&$G_{32}$\\
\hline
\end{tabular}
\end{center}
\begin{center}
\tabcolsep=0.85 mm
\begin{tabular}{|c|c|c|c|c|}
\hline
Size&Type& Cap&Weight enumerators&$G$\\
\hline
\hline
$9$&
Inc.&
{\scriptsize
\tabcolsep=0.85 mm
\begin{tabular}{ccccccccc}
 1& 1& 0& 0& 0& 0& 0& 1& 1\\
 0& 1& 1& 1& 0& 0& 0& 1& 0\\
 0& 1& 0& 1& 1& 0& 0& 0& 1\\
 0& 0& 0& 1& 0& 1& 0& 1& 1\\
 0& 0& 0& 0& 0& 0& 1& 0& 0\\
\end{tabular}
}
&$1^{1}4^{14}5^{14}8^{1}9^{1}$&$G_{1334}$\\
\hline
$9$&
Inc.&
{\scriptsize
\tabcolsep=0.85 mm
\begin{tabular}{ccccccccc}
 1& 1& 1& 0& 0& 0& 0& 1& 1\\
 0& 0& 1& 1& 0& 0& 0& 1& 0\\
 0& 1& 1& 0& 1& 0& 0& 0& 1\\
 0& 0& 0& 0& 0& 1& 0& 1& 1\\
 0& 1& 0& 0& 0& 0& 1& 1& 0\\
\end{tabular}
}
&$3^{4}4^{14}5^{8}7^{4}8^{1}$&$G_{192}$\\
\hline
$9$&
Inc.&
{\scriptsize
\tabcolsep=0.85 mm
\begin{tabular}{ccccccccc}
 1& 1& 0& 0& 0& 0& 0& 1& 0\\
 0& 1& 1& 0& 0& 0& 0& 0& 1\\
 0& 1& 0& 1& 1& 0& 0& 1& 0\\
 0& 0& 0& 0& 1& 1& 0& 1& 1\\
 0& 0& 0& 0& 1& 0& 1& 0& 1\\
\end{tabular}
}
&$3^{6}4^{9}5^{9}6^{6}9^{1}$&$G_{72}$\\
\hline
$9$&
Inc.&
{\scriptsize
\tabcolsep=0.85 mm
\begin{tabular}{ccccccccc}
 1& 1& 1& 0& 0& 0& 0& 1& 0\\
 0& 0& 1& 1& 0& 0& 0& 0& 1\\
 0& 1& 1& 0& 1& 0& 0& 1& 0\\
 0& 0& 0& 0& 0& 1& 0& 1& 1\\
 0& 1& 0& 0& 0& 0& 1& 0& 1\\
\end{tabular}
}
&$2^{1}3^{3}4^{11}5^{11}6^{3}7^{1}9^{1}$&$G_{48}$\\
\hline
$9$&
Comp.&
{\scriptsize
\tabcolsep=0.85 mm
\begin{tabular}{ccccccccc}
 1& 1& 0& 0& 0& 0& 0& 1& 1\\
 0& 1& 1& 1& 0& 0& 0& 1& 0\\
 0& 1& 0& 1& 1& 0& 0& 0& 1\\
 0& 0& 0& 1& 0& 1& 0& 1& 1\\
 0& 0& 0& 0& 0& 0& 1& 1& 0\\
\end{tabular}
}
&$2^{1}4^{21}6^{7}8^{2}$&$G_{336}$\\
\hline
\hline
$10$&
Inc.&
{\scriptsize
\tabcolsep=0.85 mm
\begin{tabular}{cccccccccc}
 1& 1& 0& 0& 0& 0& 0& 1& 1& 0\\
 0& 1& 1& 0& 0& 0& 0& 1& 0& 1\\
 0& 1& 0& 1& 1& 0& 0& 0& 1& 0\\
 0& 0& 0& 0& 1& 1& 0& 0& 1& 1\\
 0& 0& 0& 0& 1& 0& 1& 1& 0& 1\\
\end{tabular}
}
&$4^{15}6^{15}10^{1}$&$G_{720}$\\
\hline
$10$&
Inc.&
{\scriptsize
\tabcolsep=0.85 mm
\begin{tabular}{cccccccccc}
 1& 0& 0& 0& 0& 1& 1& 0& 1& 1\\
 0& 1& 0& 0& 0& 0& 0& 1& 1& 1\\
 0& 0& 1& 0& 0& 1& 1& 0& 1& 1\\
 0& 0& 0& 1& 0& 0& 1& 1& 0& 1\\
 0& 0& 0& 0& 1& 1& 0& 1& 0& 1\\
\end{tabular}
}
&$2^{1}4^{6}5^{16}6^{6}8^{1}10^{1}$&$G_{384}$\\
\hline
$10$&
Inc.&
{\scriptsize
\tabcolsep=0.85 mm
\begin{tabular}{cccccccccc}
 1& 0& 0& 0& 0& 1& 0& 1& 0& 0\\
 0& 1& 0& 0& 0& 0& 1& 1& 1& 0\\
 0& 0& 1& 0& 0& 1& 0& 1& 1& 1\\
 0& 0& 0& 1& 0& 1& 1& 0& 1& 1\\
 0& 0& 0& 0& 1& 0& 1& 0& 0& 1\\
\end{tabular}
}
&$3^{2}4^{7}5^{12}6^{7}7^{2}10^{1}$&$G_{48}$\\
\hline
$10$&
Comp.&
{\scriptsize
\tabcolsep=0.85 mm
\begin{tabular}{cccccccccc}
 1& 1& 0& 0& 0& 0& 1& 0& 1& 1\\
 1& 0& 1& 0& 0& 0& 0& 1& 0& 1\\
 0& 0& 0& 1& 0& 0& 1& 1& 1& 1\\
 1& 0& 0& 0& 1& 0& 0& 1& 1& 0\\
 1& 0& 0& 0& 0& 1& 1& 1& 0& 0\\
\end{tabular}
}
&$4^{10}5^{16}8^{5}$&$G_{1920}$\\
\hline
\hline
$11$&
Inc.&
{\scriptsize
\tabcolsep=0.85 mm
\begin{tabular}{ccccccccccc}
 1& 0& 0& 0& 0& 1& 1& 0& 1& 0& 1\\
 0& 1& 0& 0& 0& 0& 0& 1& 1& 1& 1\\
 0& 0& 1& 0& 0& 1& 1& 0& 1& 1& 1\\
 0& 0& 0& 1& 0& 0& 1& 1& 0& 1& 1\\
 0& 0& 0& 0& 1& 1& 0& 1& 0& 0& 1\\
\end{tabular}
}
&$3^{1}4^{2}5^{12}6^{12}7^{2}8^{1}11^{1}$&$G_{192}$\\
\hline
$11$&
Inc.&
{\scriptsize
\tabcolsep=0.85 mm
\begin{tabular}{ccccccccccc}
 1& 1& 0& 0& 0& 0& 1& 0& 1& 0& 0\\
 0& 1& 1& 0& 0& 0& 0& 1& 1& 1& 0\\
 0& 0& 0& 1& 0& 0& 1& 0& 1& 1& 1\\
 0& 0& 0& 0& 1& 0& 1& 1& 0& 1& 1\\
 0& 1& 0& 0& 0& 1& 0& 1& 0& 0& 1\\
\end{tabular}
}
&$4^{5}5^{10}6^{10}7^{5}11^{1}$&$G_{120}$\\
\hline
\end{tabular}
\end{center}
\begin{center}
\tabcolsep=0.85 mm
\begin{tabular}{|c|c|c|c|c|}
\hline
Size&Type& Cap&Weight enumerators&$G$\\
\hline
\hline
$12$&
Inc.&
{\scriptsize
\tabcolsep=0.85 mm
\begin{tabular}{cccccccccccc}
 1& 0& 0& 0& 0& 1& 1& 0& 1& 1& 0& 1\\
 0& 1& 0& 0& 0& 0& 0& 1& 1& 1& 1& 1\\
 0& 0& 1& 0& 0& 1& 1& 0& 0& 1& 1& 1\\
 0& 0& 0& 1& 0& 0& 1& 1& 1& 0& 1& 1\\
 0& 0& 0& 0& 1& 1& 0& 1& 0& 0& 0& 1\\
\end{tabular}
}
&$4^{3}6^{24}8^{3}12^{1}$&$G_{2304}$\\
\hline
$12$&
Inc.&
{\scriptsize
\tabcolsep=0.85 mm
\begin{tabular}{cccccccccccc}
 1& 1& 0& 0& 0& 0& 1& 0& 1& 0& 0& 1\\
 0& 1& 1& 0& 0& 0& 0& 1& 1& 1& 0& 1\\
 0& 0& 0& 1& 0& 0& 1& 0& 1& 1& 1& 1\\
 0& 0& 0& 0& 1& 0& 1& 1& 0& 1& 1& 1\\
 0& 1& 0& 0& 0& 1& 0& 1& 0& 0& 1& 1\\
\end{tabular}
}
&$4^{1}5^{8}6^{12}7^{8}8^{1}12^{1}$&$G_{192}$\\
\hline
\hline
$13$&
Inc.&
{\scriptsize
\tabcolsep=0.85 mm
\begin{tabular}{ccccccccccccc}
 1& 0& 0& 0& 0& 1& 1& 0& 1& 1& 0& 0& 1\\
 0& 1& 0& 0& 0& 0& 0& 1& 1& 1& 1& 0& 1\\
 0& 0& 1& 0& 0& 1& 1& 0& 0& 1& 1& 1& 1\\
 0& 0& 0& 1& 0& 0& 1& 1& 1& 0& 1& 1& 1\\
 0& 0& 0& 0& 1& 1& 0& 1& 0& 0& 0& 1& 1\\
\end{tabular}
}
&$5^{3}6^{12}7^{12}8^{3}13^{1}$&$G_{576}$\\
\hline
\hline
$14$&
Inc.&
{\scriptsize
\tabcolsep=0.85 mm
\begin{tabular}{cccccccccccccc}
 1& 1& 0& 0& 0& 0& 1& 1& 0& 1& 1& 0& 0& 1\\
 0& 1& 1& 0& 0& 0& 0& 0& 1& 1& 1& 1& 0& 1\\
 0& 0& 0& 1& 0& 0& 1& 1& 0& 0& 1& 1& 1& 1\\
 0& 0& 0& 0& 1& 0& 0& 1& 1& 1& 0& 1& 1& 1\\
 0& 1& 0& 0& 0& 1& 1& 0& 1& 0& 0& 0& 1& 1\\
\end{tabular}
}
&$6^{7}7^{16}8^{7}14^{1}$&$G_{2688}$\\
\hline
\hline
$15$&
Inc.&
{\scriptsize
\tabcolsep=0.85 mm
\begin{tabular}{ccccccccccccccc}
 1& 1& 0& 0& 0& 0& 1& 1& 1& 0& 1& 1& 0& 0& 1\\
 0& 1& 1& 0& 0& 0& 0& 0& 0& 1& 1& 1& 1& 0& 1\\
 0& 0& 0& 1& 0& 0& 1& 1& 0& 0& 0& 1& 1& 1& 1\\
 0& 0& 0& 0& 1& 0& 0& 1& 1& 1& 1& 0& 1& 1& 1\\
 0& 1& 0& 0& 0& 1& 1& 0& 1& 1& 0& 0& 0& 1& 1\\
\end{tabular}
}
&$7^{15}8^{15}15^{1}$&$G_{20160}$\\
\hline
\hline
$16$&
Comp.&
{\scriptsize
\tabcolsep=0.85 mm
\begin{tabular}{cccccccccccccccc}
 1& 1& 0& 0& 0& 0& 1& 1& 1& 0& 1& 1& 0& 0& 0& 1\\
 0& 1& 1& 0& 0& 0& 0& 0& 0& 1& 1& 1& 1& 1& 0& 1\\
 0& 0& 0& 1& 0& 0& 1& 1& 0& 0& 0& 1& 1& 1& 1& 1\\
 0& 0& 0& 0& 1& 0& 0& 1& 1& 1& 1& 0& 0& 1& 1& 1\\
 0& 1& 0& 0& 0& 1& 1& 0& 1& 1& 0& 0& 1& 0& 1& 1\\
\end{tabular}
}
&$8^{30}16^{1}$&$G_{322560}$\\
\hline
\end{tabular}
\end{center}


\begin{thebibliography}{7}
\small

\bibitem{18secondo} A. Davydov, G. Faina, S. Marcugini and F. Pambianco, On size of complete caps in projective spaces $PG(n,q)$ and arcs in planes $PG(2,q)$, \emph{Journal of Geometry, published online} 18/7/2009.
\bibitem{calottecomplete} G. Faina, S. Marcugini, A. Milani and F. Pambianco, The size $k$ of the complete $k$-caps in $PG(n,q)$ for small $q$ and $3 \leq n \leq 5$,  \emph{Ars Combinatoria} \textbf{50}  (1998), 235-243.
\bibitem{Hall} M. Hall and J. K. Senior, \emph{The group of order $2^{n}$ ($n\leq 6$)}, Macmillan, New York, 1964.
\bibitem{Hir} J. W. P. Hirschfeld, \emph{Projective geometries over finite fields}, Claredon Press, Oxford, 1979.
\bibitem{Hir2} J. W. P. Hirschfeld, \emph{Finite projective spaces of three dimension}, Claredon Press, Oxford, 1985.





\end{thebibliography}
\end{document}